# Linear-time Interval Algorithm For Time-varying Power Flow

**Zongjie Wang and C. Lindsay Anderson**

Department of Biological and Environmental Engineering, Cornell University, Ithaca, NY, 14850. zw337@cornell.edu (W.Z.); cla28@cornell.edu (A.C.)

**Abstract:**

Balanced steady operation state in power networks changes over time. Traditional power flow algorithm is focused on the steady operation state under a certain time point and calculates the corresponding voltage and power distributions for given nodal power injections and network topology. A linear-time interval algorithm for time-varying power flow is first proposed in this paper, which is focused on the balanced operation state under a time interval and computes the voltage and power distributions at any time point via a time-varying function, through which the power flow time-varying process is well represented in steady-date power networks. In this paper, a linear-time interval regarding nodal power injections is defined, in which the norms of nodal voltage derivatives and Jacobian matrix are analyzed, a combined time-varying nodal voltage function is proposed, and a linear-time interval algorithm for time-varying power flow with high accuracy is finally achieved, which is leveraged for simplifying the complexity of solving non-linear dynamic time-varying problems. Simulation case studies on IEEE 5-bus and modified 118-bus systems have demonstrated the effectiveness and efficiency of the proposed algorithm.

## 1. Introduction

Power flow analysis is the most important and essential approach to investigate problems in power system operation and planning [1-3]. Based on a specified generating state and transmission network structure, power flow analysis solves the steady operation state with nodal voltages and branch power flow in power systems [4-6].

Power flow analysis is an important prerequisite in power systems, such as fault analysis, stability studies, economic operation, and optimal power flow [7-10]. Power flow could be a support during examining effectiveness of the alternative plans for future network expansion, when adding new generators or transmission lines are needed [11]. It also helps in continuously monitoring of the current state of power systems, thus it is applied on daily basis in load dispatch/optimization [12, 13].

The traditional optimal power flow (TOPF) is focused on the power balancing of a specified time point, and assumes that the OPF solution at that specified time point will satisfy the system constraints until the next time point. However, since the power steady operation states in power systems are time-varying, the power flow distributions at the specified time points may not perfectly represent other non-specified time points. In addition, variable renewable energies (VREs) including wind farms and solar sites have been rapidly integrated into power systems and brings into a lot of uncertainty with more fluctuant and unpredictable loading conditions, which may jeopardize the security and stability of the power systems. As a result, there is no guarantee that solutions obtained at specific time points will satisfy constraints at non-specified time points, especially for binding constraints. An efficient OPF approach should ensure that all the load balance, thermal, and voltage limit constraints are satisfied at any time point within a given time interval. More specifically, a dynamic optimal power flow (DOPF) model was proposed in [14] as an extension of TOPF which aims to cover multiple time periods and provides a more realistic OPF model, leading to better dispatch decisions across time periods. Although DOPF decomposes the entire time period into multiple time steps, it tends to be computationally intensive since the sets of constraints over all the time points are taken into account [13,15-17].

To specify proper effective discrete time points, to ease the computational burden and to effectively represent the optimal power flow during a continuous time period, as a prerequisite of the above, this paper first proposes the concept of linear-time interval, in which the power injections are linear mapping of time. Linear-time interval is a fundamental and important concept through this research, it actually exists in power systems and represents the forecasting output level between two consecutive discrete time points during generation scheduling plans. Under a linear-time interval, properties of the norms of voltage derivatives and Jacobian matrix are analyzed; the nodal voltages are theoretically shown to be approximately linear of time, and this property even holds for heavy power variation ranges that are above 100%. In addition, a combined time-varying nodal voltage function is further proposed, in which the maximum nodal voltage errors are ensured to be restricted within $10^{-4}$. A linear-time interval algorithm for time-varying power flow with high accuracy is finally presented and evaluated, which is leveraged for simplifying the complexity of solving non-linear dynamic time varying problems.

The remainder of this paper is organized as follows. In Section 2, we review the power forecasting output and propose the linear-time interval. The derivatives of nodal voltages are given in Section 3. Some voltage properties under linear-time interval are discussed and verified in a IEEE 5-bus system in Section 4 wherein the norms of Jacobian matrix and nodal voltage derivatives are studied. In Section 5, a combined time-varying nodal voltage function is further proposed and evaluated in the same IEEE 5-bus system. In Section 6, a linear-time interval algorithm for time-varying power flow is described, simulation case studies in a modified IEEE 118-bus system are then recorded. Conclusions and future work regarding optimal power flow are finally given in Sections 9 and 10.

## 2. Linear-time intervals

We first assume the beginning time point and ending time point of a time interval $T_l$ are $t_0$, $t_e$ respectively, which is $T_l \in [t_0, t_e]$. For a time interval $T_l$, if the active power of all nodes (except the slack node) and reactive power of all PQ nodes are linear functions of time $t$, then we define $T_l$ as a linear-time interval.

The defined linear-time interval regarding power injections actually exists in power systems.

From engineering perspective, the dispatch control center arranges the generation scheduling plans based on load forecasts/equivalent load forecasts with VREs, such as wind/solar power sources. Load forecasts including total and nodal loads are oriented towards discrete time points for a future time period. A straight line segment is achieved by connecting the load forecasts between two consecutive discrete time points, thus the load forecasting output and the corresponding generation schedules for a future time period is a piecewise linear function of time. Figure 1 gives an example of nodal power forecasting output with hourly discrete forecasting time points during day-ahead scheduling, in which the two power outputs could either represent active or reactive power.

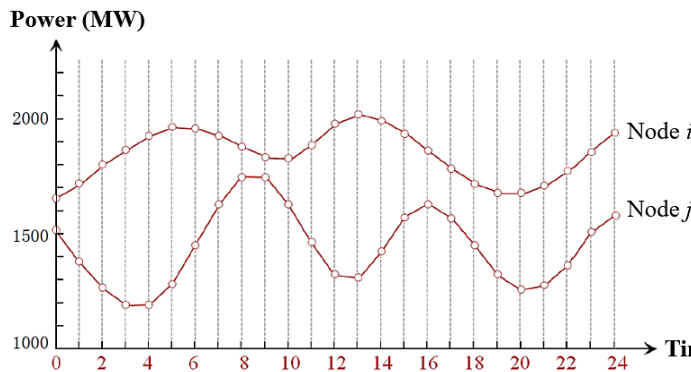

**Figure 1.** Nodal forecasting power output during day-ahead scheduling.

From mathematic perspective, all non-linear time-varying curves can be approximated as piecewise linear functions combined with many linear-time intervals, as the intervals get smaller/the number of intervals gets bigger, the piecewise linear functions gradually match closer to the non-linear time-varying curves. Therefore, linear-time interval is the fundamental approach to solve steady-state non-linear time-varying problems.

The actual load forecasting output is a time-varying curve instead of a piecewise linear function. While making generation scheduling plans for future dispatch time period, dispatch control center only has information about the piecewise linear function of load forecasting. Just like the load forecasting error is inevitable, the error caused from replacing with piecewise linear function is also inevitable.

The proposed linear-time interval discretize a continuous non-linear time-varying problem into combined linear-time intervals sub-problems, and thus reduces the complexity of finding solutions under non-linear time-varying problems.

If we assume $\Delta T_l$ is the length of the linear-time interval $T_l$, and $\Delta t$ is the time difference between any time point $t$ and the beginning time point $t_0$, then

$$\begin{cases} \Delta T_l = t_e - t_0 \\ \Delta t = t - t_0 \end{cases}. \tag{1}$$

Thus, for any time point, $t \in T_l$, the active power of all nodes (except the slack node) and reactive power of all PQ nodes can be represented as follows:

$$\begin{cases} P_i(t) = P_i(t_0) + K_{Pi} \Delta t, & (\forall i \notin V\theta \text{ bus}) \\ Q_i(t) = Q_i(t_0) + K_{Qi} \Delta t, & (\forall i \in PQ \text{ bus}) \end{cases}, \tag{2}$$

where $K_{Pi}$, $K_{Qi}$ are the slopes of active and reactive power functions with respect to $t$, respectively:

$$\begin{cases} K_{Pi} = \dfrac{1}{\Delta T_l} \Delta P_i \\ K_{Qi} = \dfrac{1}{\Delta T_l} \Delta Q_i \end{cases}. \tag{3}$$

Herein, $\Delta P_i$ and $\Delta Q_i$ are respectively the increments of active and reactive power of node $i$ for linear-time interval $T_l$:

$$\begin{cases} \Delta P_i = P_i(t_e) - P_i(t_0) \\ \Delta Q_i = Q_i(t_e) - Q_i(t_0) \end{cases}. \tag{4}$$

The vector form of (2) is expressed as

$$\begin{cases} \boldsymbol{P}(t) = \boldsymbol{P}(t_0) + \Delta t \boldsymbol{K}_P \\ \boldsymbol{Q}(t) = \boldsymbol{Q}(t_0) + \Delta t \boldsymbol{K}_Q \end{cases}. \tag{5}$$

If

$$\boldsymbol{y}(t) = \begin{bmatrix} \boldsymbol{P}(t) \\ \boldsymbol{Q}(t) \end{bmatrix}, \boldsymbol{k} = \begin{bmatrix} \boldsymbol{K}_P \\ \boldsymbol{K}_Q \end{bmatrix}, \tag{6}$$

then

$$\boldsymbol{y}(t) = \boldsymbol{y}(t_0) + \Delta t \boldsymbol{k}, \tag{7}$$

where $\boldsymbol{k}$ is the slope vector for linear-time interval $T_l$ and is related to different time units. In this paper, the slope vectors in hour, minute and second units are represented as $\boldsymbol{k}_h, \boldsymbol{k}_m, \boldsymbol{k}_s$ and satisfy the following condition:

$$\boldsymbol{k}_h = 60 \boldsymbol{k}_m = 3600 \boldsymbol{k}_s. \tag{8}$$

Without loss of generality, we choose hourly unit in this paper with the corresponding slope vector $k_h$.

## 3. Nodal voltage derivatives

In rectangular coordinate, the time-varying nodal voltage vector $x(t)$ is

$$x(t) = \begin{bmatrix} V_{re}(t) \\ V_{im}(t) \end{bmatrix}, \quad (9)$$

where $V_{re}(t)$ and $V_{im}(t)$ are the real part and imaginary part of nodal voltage vectors, respectively. Thus the time-varying power equation is expressed as:

$$h(x(t)) = y(t). \quad (10)$$

The corresponding first derivative of nodal voltages is

$$J(t)x^{(1)}(t) = k, \quad (11)$$

where $J(t)$ is the Jacobian matrix shown as follows:

$$J(t) = \frac{\partial}{\partial x} h(x(t)). \quad (12)$$

Since the elements in Jacobian matrix $J(t)$ are linear function of nodal voltages in rectangular coordinates, we get

$$J^{(k)}(t) = J(x^{(k)}(t)), (k \geq 0). \quad (13)$$

By continuously calculating the derivatives of (11), we get the following $d^{th}$ order of nodal voltage derivative:

$$J(t)x^{(d)}(t) = b_d(t), (d \geq 1), \quad (14)$$

herein, when $d = 1$:

$$b_1(t) = k; \quad (15)$$

when $d \geq 2$:

$$b_d(t) = -\sum_{k=1}^{d-1} C_{d-1}^k J(x^{(k)}(t))x^{(d-k)}(t), (d \geq 2), \quad (16)$$

where

$$C_{d-1}^k = \frac{(d-1)!}{k!(d-1-k)!}. \quad (17)$$

## 4. Time-varying nodal voltage properties in linear-time intervals

*4.1 Matrix and vector norms*

Norms are compatible, if

$$mv = b, \quad (18)$$

then

$$\|mv\|_p = \|b\|_p, \quad (19)$$

and satisfies

$$\|m\|_p \|v\|_p \geq \|b\|_p, \quad (20)$$

where $m$ is matrix, $v$ and $b$ are the vectors that are compatible with matrix $m$.

Norms also have equivalence property, for $p$-norms when $p = 1, 2, \infty$, if one of the norms follows the corresponding equality/inequality, so do the other two norms.

*4.2 Norm of the derivatives of Jacobian matrix*

**Proposition 1.** *In rectangular coordinate, the derivative matrix of Jacobian, $J^{(k)}(t)$, satisfies the following norm relationship:*

$$\|J^{(k)}(t)\| \leq \|J(t)\|\|x^{(k)}(t)\|, (k \geq 0). \tag{21}$$

Herein, $\|\cdot\|$ could be any of the three norms for $p = 1, 2, \infty$.

**Proof.** The element in the derivative matrix of Jacobian, $J^{(k)}(t)$, is expressed as

$$J_{ij}^{(k)}(t) = \frac{\partial J_{ij}(t)}{\partial x_l} x_l^{(k)}(t), (k \geq 0). \tag{22}$$

According to the definition of norm, the $\infty$-norm is shown as follows

$$\|J^{(k)}(t)\|_\infty = \max_j \sum_j |J_{ij}^{(k)}(t)|, \tag{23}$$

obviously,

$$\|J^{(k)}(t)\|_\infty \leq \max_j \sum_j \left|\frac{\partial J_{ij}(t)}{\partial x_l}\right| \max_l |x_l^{(k)}(t)|. \tag{24}$$

Since

$$\|x^{(k)}(t)\|_\infty = \max_l |x_l^{(k)}(t)|, \tag{25}$$

we get that

$$\|J^{(k)}(t)\|_\infty \leq \|x^{(k)}(t)\|_\infty \max_j \sum_j \left|\frac{\partial J_{ij}(t)}{\partial x_l}\right|. \tag{26}$$

Note that

$$\|J(t)\|_\infty \approx \max_j \sum_j \left|\frac{\partial J_{ij}(t)}{\partial x_l}\right|, \tag{27}$$

therefore

$$\|J^{(k)}(t)\|_\infty \leq \|J(t)\|_\infty \|x^{(k)}(t)\|_\infty, (k \geq 0). \tag{28}$$

According to the norm equivalence property, the above proposition also holds 1-norm and 2-norm. This ends the proof.

*4.3 Norm of nodal voltage derivative vectors*

The condition number of Jacobian matrix $J(t)$ is

$$\rho_t = \|J^{-1}(t)\|_p \|J(t)\|_p, (p = 1, 2, \infty). \tag{29}$$

Thus, the norm of Jacobian inverse matrix is expressed as

$$\|J^{-1}(t)\| = \rho_t \|J(t)\|^{-1}. \tag{30}$$

**Proposition 2.** *For a $d^{th}$-order derivative of nodal voltage vectors, $x^{(d)}(t)$, the corresponding norm satisfies the following condition:*

$$\|x^{(d)}(t)\| \leq (2d-3)!! \rho_t^{d-1} \|x^{(1)}(t)\|^d, (d \geq 2), \tag{31}$$

*where $(\cdot)!!$ represents double factorial.*

**Proof.** The mathematical induction is employed with the following equation achieved from (14) and (16):

$$x^{(2)}(t) = -J^{-1}(t)J(x^{(1)}(t))x^{(1)}(t). \tag{32}$$

According to proposition 1, we have

$$\|x^{(2)}(t)\| \leq \|J^{-1}(t)\|\|J(t)\|\|x^{(1)}(t)\|^2. \tag{33}$$

Considering (30), we get that

$$\|x^{(2)}(t)\| \leq \rho_t \|x^{(1)}(t)\|^2. \tag{34}$$

Thus the proposition holds when the order $d$ is 2.

According to (16), since

$$\boldsymbol{b}_3(t) = -2\boldsymbol{J}(\boldsymbol{x}^{(2)}(t))\boldsymbol{x}^{(1)}(t) - \boldsymbol{J}(\boldsymbol{x}^{(1)}(t))\boldsymbol{x}^{(2)}(t), \tag{35}$$

then

$$\|\boldsymbol{b}_3(t)\| \leq 3\|\boldsymbol{J}(t)\|\|\boldsymbol{x}^{(2)}(t)\|\|\boldsymbol{x}^{(1)}(t)\|, \tag{36}$$

Considering (34), we get that

$$\|\boldsymbol{b}_3(t)\| \leq 3\rho_t \|\boldsymbol{J}(t)\|\|\boldsymbol{x}^{(1)}(t)\|^3. \tag{37}$$

Obviously, the following condition is satisfied

$$\|\boldsymbol{x}^{(3)}(t)\| \leq 3\rho_t^2 \|\boldsymbol{x}^{(1)}(t)\|^3. \tag{38}$$

Thus the proposition also holds when the order $d$ is 3.

We assume that the proposition holds in $(d-1)^{th}$-order of derivative. When the order is $d$, we obtain the following condition from (16):

$$\|\boldsymbol{b}_d(t)\| \leq \sum_{k=1}^{d-1} C_{d-1}^k \|\boldsymbol{J}(\boldsymbol{x}^{(d-k)}(t))\|\|\boldsymbol{x}^{(k)}(t)\|, \tag{39}$$

According to proposition 1, we get that

$$\|\boldsymbol{b}_d(t)\| \leq \sum_{k=1}^{d-1} C_{d-1}^k \|\boldsymbol{J}(t)\|\|\boldsymbol{x}^{(d-k)}(t)\|\|\boldsymbol{x}^{(k)}(t)\|. \tag{40}$$

Since the proposition holds in $(d-1)^{th}$-order, the following two conditions are satisfied:

$$\|\boldsymbol{x}^{(d-k)}(t)\| \leq (2(d-k)-3)!!\,\rho_t^{d-k-1} \|\boldsymbol{x}^{(1)}(t)\|^{d-k}, \tag{41}$$

$$\|\boldsymbol{x}^{(k)}(t)\| \leq (2k-3)!!\,\rho_t^{k-1} \|\boldsymbol{x}^{(1)}(t)\|^k. \tag{42}$$

Hence

$$\|\boldsymbol{b}_d(t)\| \leq \varphi \rho_t^{d-2} \|\boldsymbol{J}(t)\|\|\boldsymbol{x}^{(1)}(t)\|^d, \tag{43}$$

herein

$$\varphi = \sum_{k=1}^{d-1} C_{d-1}^k (2k-3)!!(2(d-k)-3)!!. \tag{44}$$

According to the proposition 3 (double factorial summation formula) in Appendix, we have

$$\varphi = (2d-3)!!. \tag{45}$$

Thus,

$$\|\boldsymbol{b}_d(t)\| \leq (2d-3)!!\,\rho_t^{d-2} \|\boldsymbol{J}(t)\|\|\boldsymbol{x}^{(1)}(t)\|^d. \tag{46}$$

(31) is then achieved by substituting (46) into the following condition

$$\|\boldsymbol{x}^{(d)}(t)\| \leq \|\boldsymbol{J}^{-1}(t)\|\|\boldsymbol{b}_d(t)\|. \tag{47}$$

This ends the proof.

The condition number $\rho_t$ can be corrected for a given proper coefficient $\alpha_t \in (0,1)$:

$$\tilde{\rho}_t = \alpha_t \rho_t \leq \rho_t, \tag{48}$$

thus the inequality condition from (31) in the proposition 2 is further transformed into the following approximate equality:

$$\|\boldsymbol{x}^{(d)}(t)\| \approx (2d-3)!!\,\tilde{\rho}_t^{d-1} \|\boldsymbol{x}^{(1)}(t)\|^d, \ (d \geq 2). \tag{49}$$

*4.4 Time-varying nodal voltage properties*

**Property 1.** *For a linear-time interval, the norms of nodal voltage derivatives present a "U" curve as the order of derivative increases. In other words, the norms first decrease rapidly then decrease slightly; followed by a slight increase then a rapid increase.*

The corresponding discussions are as follows.

In the per-unit system, the order of magnitude of slope vector norm $\|k_h\|$ is within $10^0$, the order of magnitude of $\|J^{-1}(t)\|$ is generally $10^{-1}$. Since

$$\|x^{(1)}(t)\| \leq \|J^{-1}(t)\|\|k_h\|, \tag{50}$$

thus the order of magnitude of $\|x^{(1)}(t)\|$ is within $10^{-1}$.

According to (31), we have

$$\gamma_d = \frac{\|x^{(d)}(t)\|}{\|x^{(d-1)}(t)\|} \approx (2d-3)\tilde{\rho}_t \|x^{(1)}(t)\|, (d \geq 2). \tag{51}$$

For a certain time point $t$, $\tilde{\rho}_t \|x^{(1)}(t)\|$ is known, thus the ratio $\gamma_d$ between adjacent orders of derivative is a linear function with respect to order $d$, and is also known as a monotonic increasing function.

The order of magnitude for $\|x^{(1)}(t)\|$ is within $10^{-1}$, since the order of magnitude for $\tilde{\rho}_t$ is $10^0$, there must exist a critical order $d_{cr}$, such that, when $d < d_{cr}$, as the norm of nodal voltage derivative gradually decreases, the decreasing trend gets slower; when $d > d_{cr}$, as the norm of nodal voltage derivative gradually increases, the decreasing trend gets faster. As a result, the norms of nodal voltage derivatives present a "U" curve as the order of nodal voltage derivative increases. In the next subsection 4.5, the curve of the norm of nodal voltage derivatives is marked as red in Figure 4.

**Property 2.** *For a linear-time interval, the nodal voltage function with respect to time has approximate linearity.*

The corresponding discussions are as follows.

According to the Taylor series expansion, the following time-varying nodal voltages in a linear-time interval are expressed as

$$x(t) = x(t_0) + \sum_{k=1} \frac{1}{k!} \Delta^k t x^{(k)}(t_0), \tag{52}$$

According to (49) from the proposition 2, we have

$$\begin{cases} \|x^{(2)}(t)\| \approx \tilde{\rho}_t \|x^{(1)}(t)\|^2 \\ \|x^{(3)}(t)\| \approx 3\tilde{\rho}_t^2 \|x^{(1)}(t)\|^3 \end{cases}. \tag{53}$$

The order of magnitude for $\|x^{(1)}(t)\|$ is within $10^{-1}$, the order of magnitude for $\tilde{\rho}_t$ is $10^0$, thus the order of magnitude for $\|x^{(2)}(t)\|$ is within $10^{-2}$; and it will be much smaller for $\|x^{(3)}(t)\|$. The descending trend could also be observed through the *"U"* curve of property 1. Under the condition of $\Delta t < 1$, if the 2$^{nd}$ and the above orders of derivative for (52) is approximated to be zero, then the nodal voltages are approximated as linear functions with respect to time, that is,

$$x(t) = x(t_0) + \Delta t x^{(1)}(t_0). \tag{54}$$

(54) is called linear time-varying function with the corresponding norm of absolute error given by

$$\|R_1(t)\| = \frac{1}{2}\Delta^2 t \|x^{(1)}(t_0 + \eta \Delta t)\|, (0 < \eta < 1). \tag{55}$$

The above absolute error is characterized through the *2$^{nd}$*-order of derivative, the corresponding order of magnitude is within $10^{-2}$. Since the orders of magnitude for real part and imaginary part of nodal voltages are $10^0$, the order of magnitude for relative voltage error is also within $10^{-2}$. This indicates that the nodal voltage function with respect to time has approximate linearity in a linear-time interval $T_l$.

*4.5 Numerical simulations of time-varying nodal voltage properties*

Based on Matpower, the IEEE 5-bus system is simulated to numerically evaluate the time-varying nodal voltage properties in a linear-time interval $T_l$. The corresponding network topology is shown in Figure 2, wherein bus 5 is a slack bus, bus 1 is a *PV* bus with a voltage magnitude of 1.05 p.u., the rest of the buses are *PQ* buses.

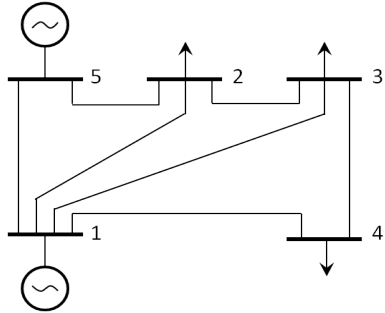

**Figure 2.** IEEE 5-bus test system.

The linear-time interval is considered as 1h. The active power output from *PQ* buses and PV bus are shown in Figure 3(a). Similarly, the reactive power output from *PQ* buses are shown in Figure 3(b).

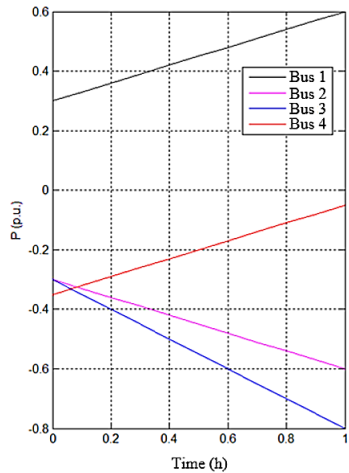 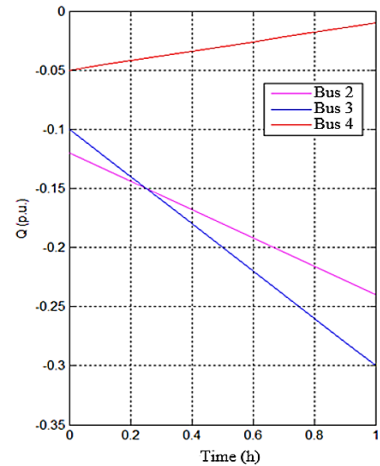

**Figure 3(a).** Active power output.    **Figure 3(b).** Reactive power output.

**Figure 3.** Power Output within a given linear-time interval.

To better illustrate the effectiveness of the proposed properties, we intentionally design some heavy node power variation ranges so that the total power variation range is above 100%. In addition, to avoid the same trends from nodal power/total power (e.g., either simultaneously increasing or decreasing), we also randomly design different trends of nodal active and reactive power, respectively (e.g., if the active power on node *i* increases, then the reactive power on node *i* decreases).

Figure 4 shows the numerical results for norms of different orders of nodal voltage derivatives in a linear-time interval $T_l$. The black curve represents the statistical values of norms of nodal voltage derivatives at different time points, which indicates that in a linear-time interval, the norms of nodal voltage derivatives present a "*U*" curve with respect to orders. The red curve represents the computational values from (49), in which the 2-norm is selected. Obviously, the values from (49) match well with the actual values and better characterize the norms of nodal voltage derivative properties.

The simulation results of time-varying nodal voltage function in a linear-time interval $T_l$ are shown in Figure 5. The "**o**" marker represents the actual nodal voltage values calculated from power flow in Matpower, the straight line between the two consecutive time points represents the computational values from linear time-varying function in (54). We can see that, the corresponding curves of real part and imaginary part of nodal voltages in a linear-time interval are extremely close to a straight line, thus can be approximated as linearity.

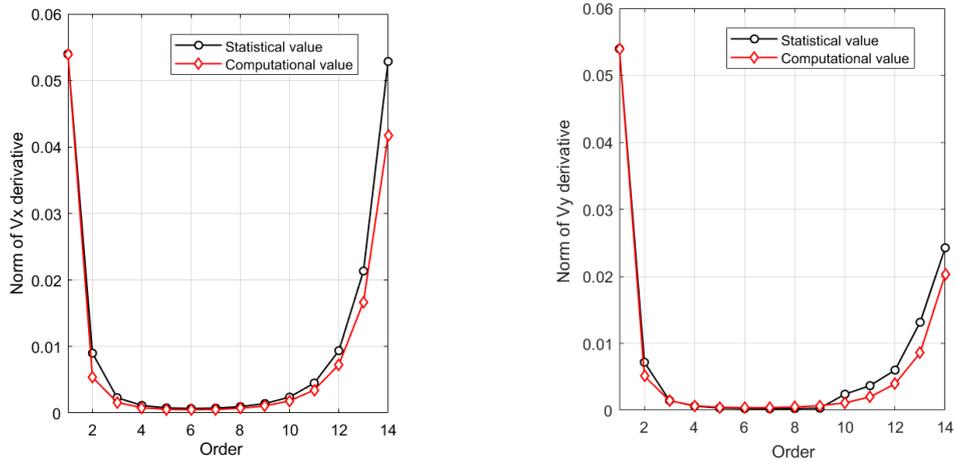

**Figure 4.** Norms of nodal voltage derivative property.

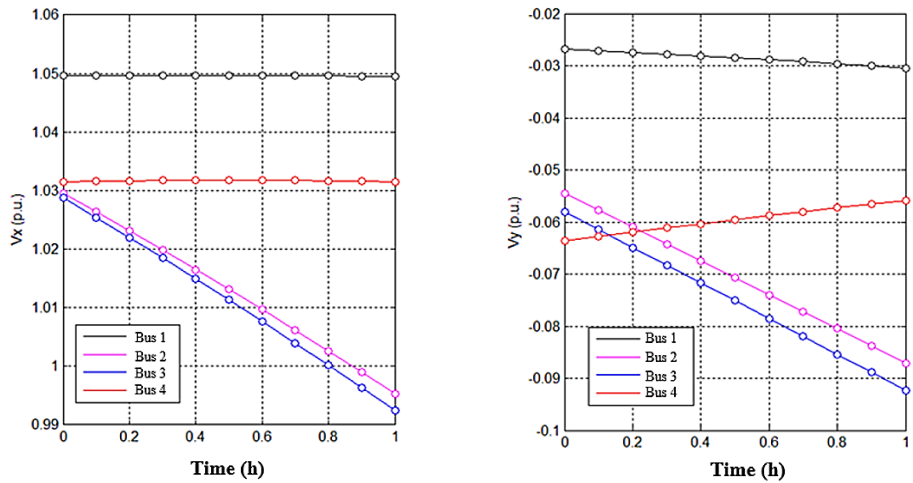

**Figure 5.** Computational results of linear time-varying function.

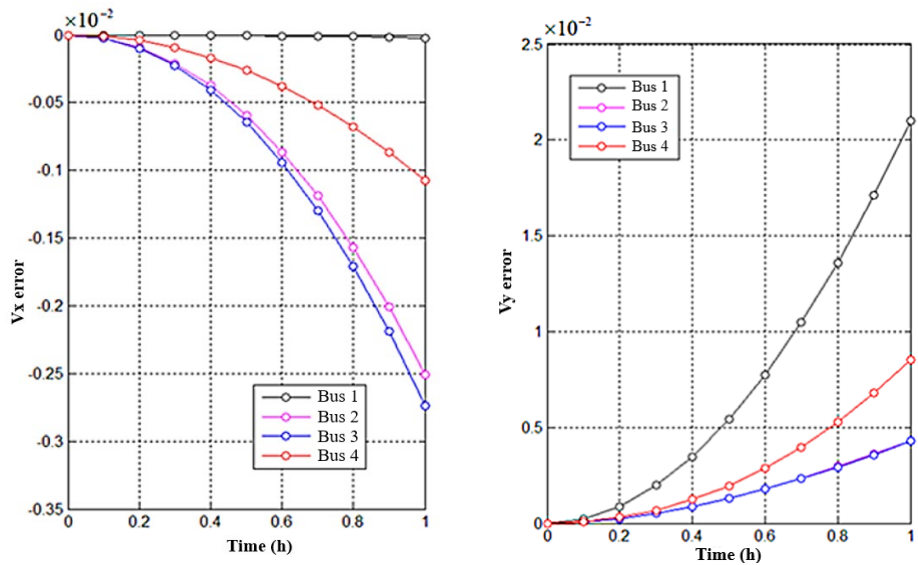

**Figure 6.** Relative errors of linear time-varying function.

Relative errors between the computational results from linear time-varying function (shown in Figure 5) and the corresponding power flow solutions calculated in Matpower (also shown in Figure

5) are further presented in Figure 6. We can see that the orders of magnitude for both real part and imaginary part of nodal voltages errors are within $10^{-2}$, which coincide with the theoretical analysis achieved from (55). Simulation results demonstrate that under a linear-time interval $T_l$, nodal voltage function is satisfied to be approximated as a linear time-varying function, as shown in (54).

## 5. Combined time-varying nodal voltage function in linear-time intervals

*5.1 Combined time-varying function*

To better evaluate the nodal voltages at any time point $t$ in linear-time interval $T_l$, the Taylor series on the beginning time point $t_0$ and the ending time point $t_e$ are respectively shown as

$$\begin{cases} x(t,t_0) = x(t_0) + \Delta t_0 x^{(1)}(t_0) + \frac{1}{2}\Delta^2 t_0 \overline{x}^{(2)} \\ x(t,t_e) = x(t_e) + \Delta t_e x^{(1)}(t_e) + \frac{1}{2}\Delta^2 t_e \overline{x}^{(2)} \end{cases}, \tag{56}$$

where

$$\begin{cases} \Delta t_0 = t - t_0 \\ \Delta t_e = t - t_e \end{cases}, \tag{57}$$

and $\overline{x}^{(2)}$ is given by the following numerical differentiation

$$\overline{x}^{(2)} = \frac{1}{\Delta T_l}(x^{(1)}(t_e) - x^{(1)}(t_0)). \tag{58}$$

(58) is the 2nd-order of derivative in (56) with the corresponding absolute error norm given by

$$\|R_2(t)\| = \frac{1}{6}\Delta^3 t \|x^{(3)}(t_0 + \eta \Delta t)\|, (0 < \eta < 1), \tag{59}$$

we can see that the order of magnitude for $\|R_2(t)\|$ is further reduced to $10^{-3}$.

Let

$$\alpha(t) = \frac{1}{\Delta T_l}\Delta t_0, \tag{60}$$

Equations in (56) is thus combined as

$$x(t) = (1-\alpha(t))x(t,t_0) + \alpha(t)x(t,t_e), \tag{61}$$

that is,

$$x(t) = x(t,t_0) + \alpha(t)(x(t,t_e) - x(t,t_0)). \tag{62}$$

Such a function is defined as a *combined time-varying nodal voltage function*.

Obviously, we can see that the difference term $x(t,t_e) - x(t,t_0)$ relates to the beginning and ending time points of a linear-time interval is evaluated in (62), which contributes to further reducing the computational errors of nodal voltage norms. Since the order of magnitude for $\alpha(t)$ is $10^{-1}$, the maximum error for (62) is thus guaranteed to be within $10^{-4}$. Compared with (54), the new proposed combined time-varying nodal voltage function is correlated to both beginning and ending time points, and is thus with a higher computational accuracy.

*5.2 Numerical simulations of combined time-varying nodal voltage function*

The same IEEE 5-bus system is tested in this subsection, followed by the same computational environments presented in subsection 4.5.

We apply (62) to compute the nodal voltages at any time point $t$ in a linear-time interval $T_l$. In Figure 7, the computational results of real part and imaginary part of nodal voltages at 11 equidistant time points in an one-hour linear-time interval are recorded. Figure 8 shows the absolute voltage errors between the computational results (shown in Figure 7) and the corresponding power flow solutions calculated in Matpower (also shown in Figure 7), through which the orders of magnitude for voltage errors are all within $10^{-6}$.

Compared with the voltage errors of the previously proposed linear time-varying function in (54) (results are shown in Figure 6), the errors shown in Figure 8 are much smaller and thus the proposed combined time-varying nodal voltage function is further illustrated and demonstrated with a higher accuracy.

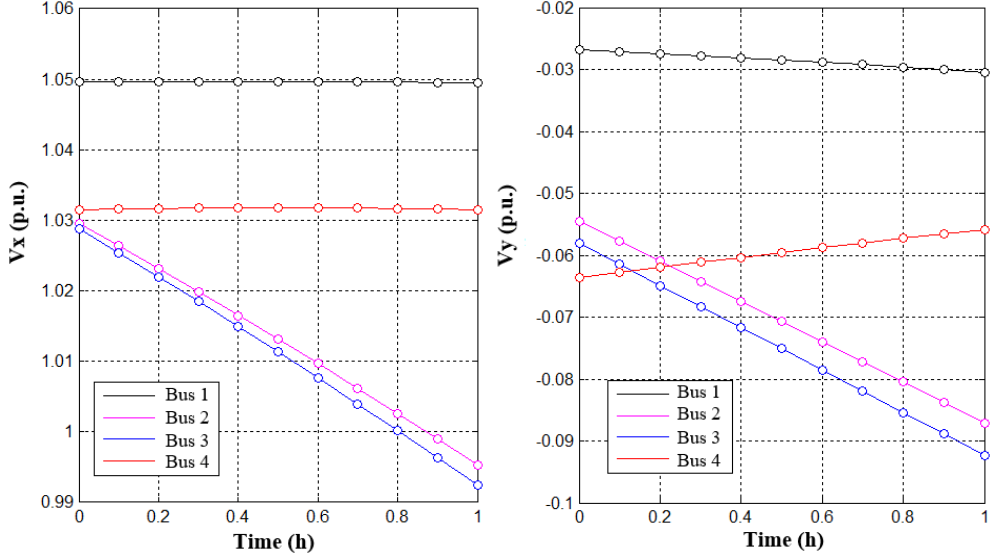

**Figure 7.** Computational results of combined time-varying nodal voltage function.

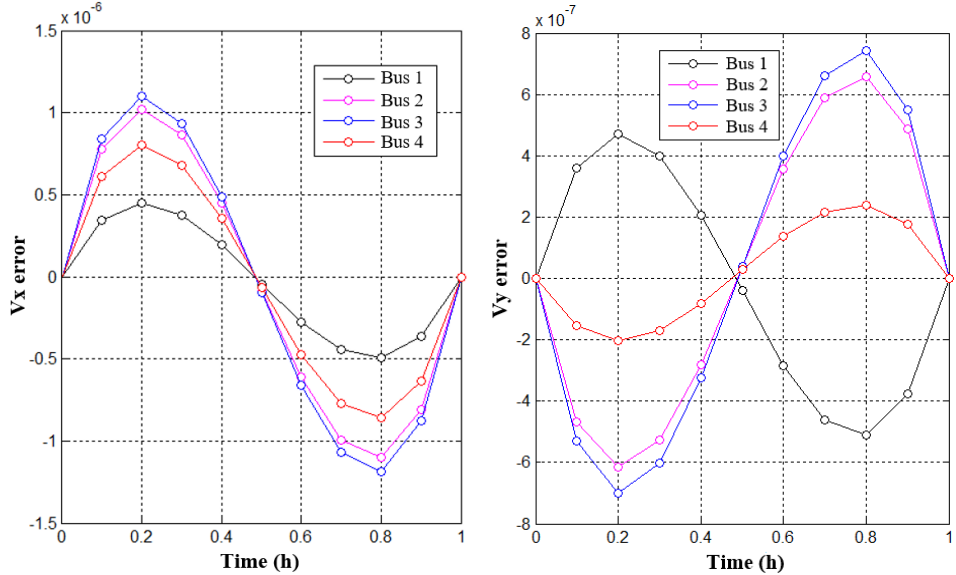

**Figure 8.** Absolute errors of combined time-varying nodal voltage function.

## 6. Linear-time interval algorithm for time-varying power flow

*6.1 Algorithm description*

Based on the proposed combined time-varying nodal voltage function, which is to calculate the time-varying functions of nodal voltages and power in a specified time period in the power systems. The corresponding steps are as follows:

*1) Linear-time interval partition*

M+1 discrete time points $t_l, (l = 0,1,...,m)$ partition the specified time period into $m$ combined linear-time intervals $T_l, (l = 1,2,...,m)$. According to (7), the piecewise linear function of nodal active and reactive power with respect to time in each linear-time interval is then achieved.

*2) Nodal voltage calculation at discrete time points*

Let the real part and imaginary part values of initial nodal voltages be 0 and 1, respectively, we first run the power flow at the beginning time point $t_0$. The initial nodal voltage values at time point $t_1$ are then obtained by linear time-varying function in (54), thus the nodal voltages at time point $t_1$ are calculated. We use the similar approach to accordingly calculate nodal voltages at other discrete time point $t_l, (l = 2,3,...,m)$ until the ending time point $t_m$.

*3) Nodal voltage calculation at linear-time intervals*

In each linear-time interval $T_l, (l = 1,2,...,m)$, we apply the combined time-varying nodal voltage function in (62) that relates to the beginning and ending time points to calculate nodal voltages in each linear-time interval $T_l$.

*4) Analysis for other parameters*

According to the combined time-varying nodal voltage function, distributions for other parameters such as branch voltages, branch power, and branch current at any time point of a given time interval can also be analyzed.

*6.2 Algorithm evaluation*

*1) Numerical stability*

Power flow analysis result at each discrete time point $t_l, (l = 2,3,...,m)$ is accurate, the nodal voltage calculations in each different linear-time interval $T_l, (l = 1,2,...,m)$ is independent, thus there does not exist the problem where the computational errors are accumulated over time. Therefore, the proposed algorithm is ideally numerically stable.

*2) Computational accuracy*

Since the norm property of different orders of nodal voltage derivatives with respect to time is given (the "U" curve), $\Delta t$ is also given, thus in each linear-time interval, the maximum errors of combined time-varying nodal voltage function are already known beforehand and ensures to be within $10^{-4}$, thus the proposed algorithm is also with high computational accuracy.

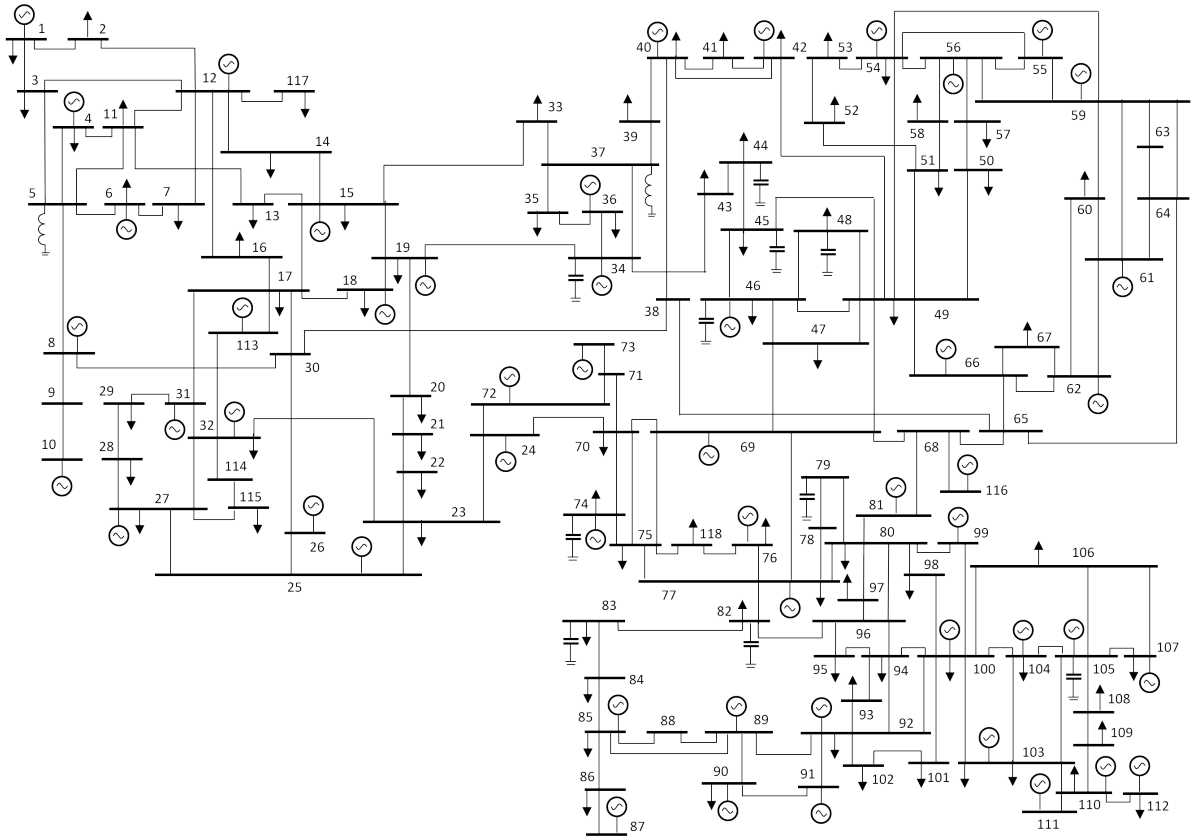

**Figure 9.** IEEE 118-bus system.

*6.3 Numerical simulations*

In this subsection, a modified IEEE 118-bus system is simulated, where 20 wind generators and 10 solar generators are added. Assume that there is no curtailment of renewable generation, then the total power generation from VREs is 30% of the total demand, with the wind power generation of being 60% of the total VREs power output. Thus a power system with high penetrations of VREs is presented. The corresponding network topology is shown in Figure 9, where bus 69 is a slack bus.

The total time period is considered as 24h, which is designed to coherent with the future dispatch period while making day-ahead generation scheduling plans. 25 discrete time points $t_l$ partition the entire time period into 24 linear-time intervals. In addition, to better evaluate and demonstrate the proposed properties, each linear-time interval (1h) is further partitioned with 11 equidistant discrete time points.

Similar to the simulations in subsection 4.5, to better illustrate the effectiveness of the proposed properties, we intentionally design some heavy nodal power variation ranges so that the total power variation range is above 100%. In addition, to avoid the same trends from nodal power/total power (e.g., either simultaneously increasing or decreasing), we also randomly design different trends of nodal active and reactive power, respectively (e.g., if the active power on node *i* increases, then the reactive power on node *i* decreases). The total active and reactive power curves are shown in Figure 10.

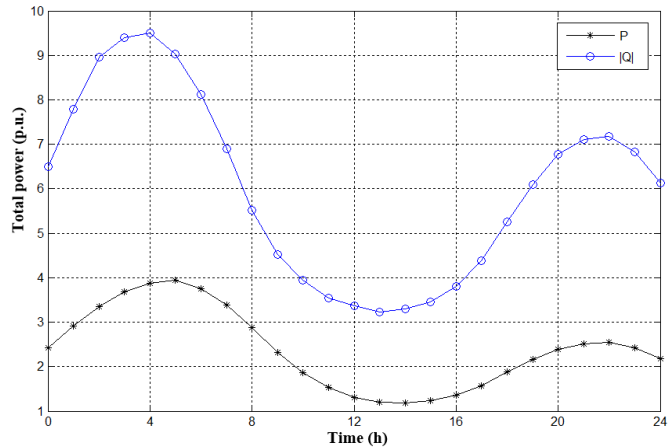

**Figure 10.** Total power curves within 24h.

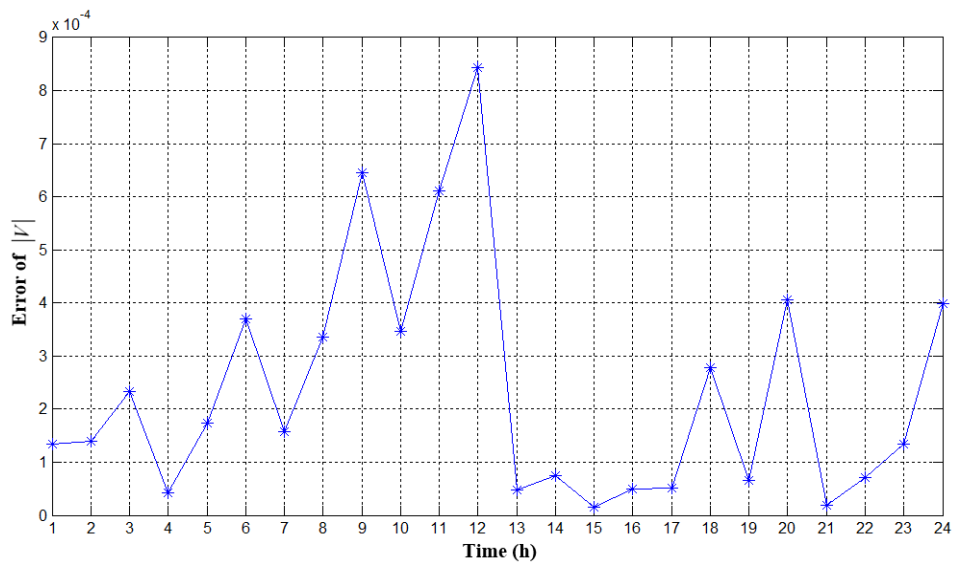

**Figure 11.** Maximum voltage errors for each linear-time interval within 24h.

Nodal voltage magnitudes at any time point in these 24 linear-time intervals are computed by (62), in other words, take node *i* for example, since there are 11 equidistant time points in each linear-time interval, the voltage magnitudes on node *i* are calculated 241 times in the 24h time period. Thus the voltage errors on each node between computational values and the corresponding power flow solutions (calculated in Matpower) at any time point in these 24 linear-time intervals are then achieved. We respectively record the maximum voltage error under each linear-time interval, with the results shown in Figure 11. Note that the maximum voltage error under each linear-time interval is selected among voltage errors of all the nodes from 11 equidistant time points, that is $118 \times 11 = 1298$.

Figure 11 shows that, the maximum nodal voltage error of all the maximum voltage errors through 24 linear-time intervals occurs on the 12th linear-time interval over the entire time period, with the corresponding error being $8.4 \times 10^{-4}$. The minimum nodal voltage error of the maximum voltage errors through 24 linear-time intervals occurs on the 15th linear-time interval over the entire time period, with the corresponding error being $1.4 \times 10^{-5}$.

We can see from the above results that, since the maximum voltage error selected among all the nodes in 11 equidistant time points in all the 24 linear-time intervals ($118 \times 241 = 28438$) is $8.4 \times 10^{-4}$, the order of magnitudes for all the nodes are guaranteed to be within $10^{-4}$, which further demonstrates the proposed linear-time interval algorithm for time-varying power flow is with high computational accuracy.

To further numerically demonstrate the effectiveness of the proposed algorithm, we next select the 21st linear-time interval for a more detailed analysis.

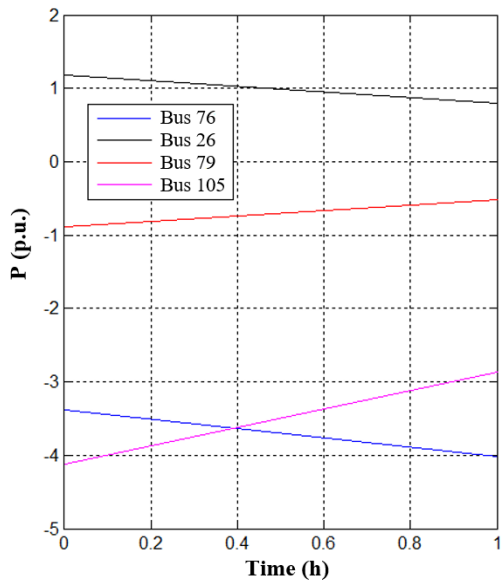 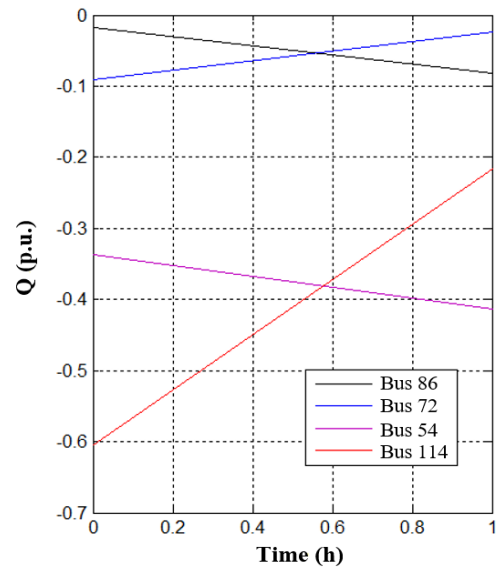

**Figure 12(a).** Active power line segments.   **Figure 12(b).** Reactive power line segments.

**Figure 12.** Load power line segments.

For the 21st linear-time interval $T_l$, we record the first four nodes with the heaviest active power variation ranges, the corresponding active power is shown in Figure 12(a). Similarly, the first four nodes with the heaviest reactive power variation ranges in $T_l$ are selected with the corresponding reactive power shown in Figure 12(b).

According to (62), the nodal voltage magnitudes at any time point of the 21st linear-time interval $T_l$ are calculated. We record four nodes with the first four maximum voltage computational errors, which are nodes 9, 16, 10 and 117. Based on the four nodes, we then record the corresponding computational values of real part and imaginary part of voltage magnitudes as shown in Figure 13. The corresponding voltage errors of real part and imaginary part are finally carried out in Figure 14.

In Figure 13, the "**o**" markers represent the actual values of nodal voltages obtained from power flow in Matpower, the straight line between two consecutive discrete time points represents the computational values from the proposed combined time-varying nodal voltage function in (62). We can see from the results that, the corresponding curves of real part and imaginary part for the selected four node voltages through the 24h time period are extremely close to straight lines, thus have approximated linearity.

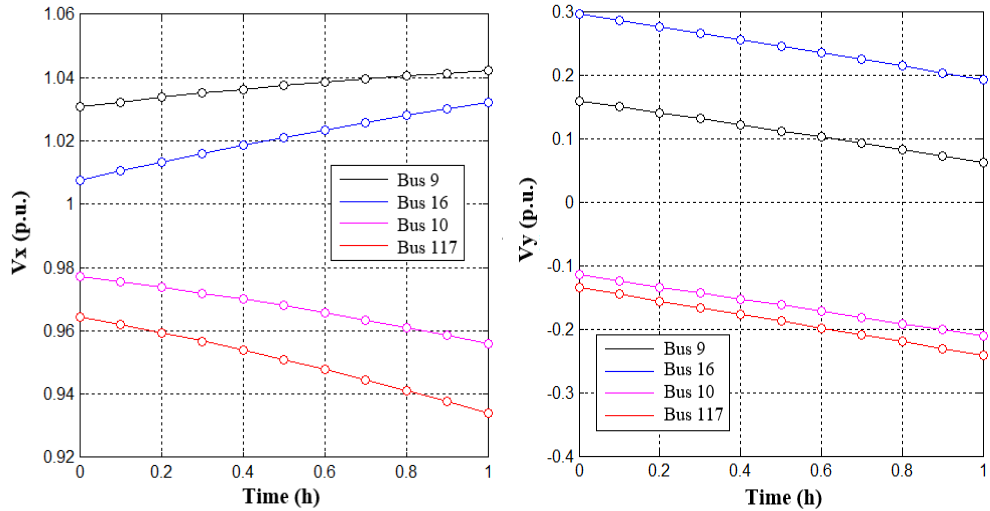

**Figure 13.** Voltage computational results from combined time-varying node voltage functions.

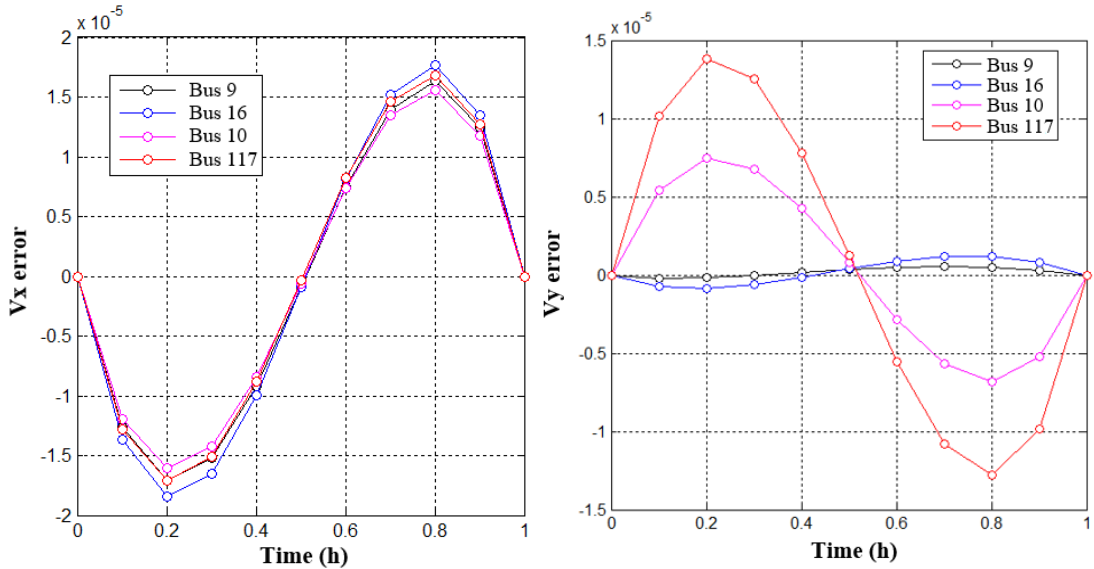

**Figure 14.** First four maximum voltage errors from combined time-varying node voltage function.

We can see from Figure 14 that, since the orders of magnitude for real-part and imaginary-part of the first four maximum voltage errors in the 21st linear-time interval $T_l$ are all within $10^{-5}$, which means all the nodal voltage errors at any time point of $T_l$ are ensured to be within $10^{-5}$, these results coincide with the theoretical analysis. Therefore, in a linear-time interval $T_l$, the nodal voltages can be approximated as a combined time-varying node voltage function with respect to time.

## 7. Conclusions

Power flow analysis is an important prerequisite in power systems especially when it comes to economic dispatch/optimal power flow. Traditional optimal power flow describes the system performance only in a single specified time point while the resulting decisions are applied to an entire time period. Since power system balanced steady operation state varies over time, it is necessary to explore the properties regarding state/control variables to ensure system's stability and security. Starting from this point of view, the conclusions of this paper are shown as a prerequisite of the above.

(i) A linear-time interval in which the nodal power injections are linear functions of time is proposed in this paper. Theoretically, linear-time interval is the fundamental approach to solve steady-state non-linear time-varying problems. Practically, linear-time interval exists in real power systems to determine the generation scheduling plans based on load forecasts. Thus linear-time interval is the key to discretize and linearize the non-linear time-varying problems.

(ii) For a linear-time interval, the norms of nodal voltage derivatives present a "U" curve with respect to the orders.

(iii) For a linear-time interval, nodal voltage magnitudes are approximate linear functions of time.

(iv) Based on the nodal voltages on the beginning and ending time points of a linear-time interval, a combined time-varying nodal voltage function is further proposed. It decouples the voltage related connections in different nodes, as well as the related connections between real part and imaginary part of nodal voltages. Namely, the proposed combined time-varying nodal voltage function easily calculates the real part and imaginary part of any node voltage magnitude at any time point of a given linear time interval with the corresponding voltage errors to be guaranteed within $10^{-4}$, thus has high promising computational accuracy.

(v) Based on the combined time-varying nodal voltage function, a linear-time interval algorithm for time-varying power flow is finally proposed. The proposed algorithm significantly simplifies the complexity of finding solutions under non-linear time-varying equations, thus is an effective, reliable and fast-responsive algorithm.

## 8. Future work

Based on the time-varying nodal voltage properties and the proposed algorithm in this paper, our future work will be focused on optimal power flow towards a time period and the corresponding applications into economic dispatch systems including day-ahead scheduling and real-time dispatch, typically in large-scale actual/synthetic power systems with high penetrations of VREs.

**Appendix**

This Appendix contains proofs of the proposition 2 given in the text above.

**Proposition 3. (Double factorial summation formula)** *For any positive integer $i \geq 2$, the following double factorial summation formula holds:*

$$(2i-1)!! = \sum_{k=1}^{i} C_i^k (2k-3)!!(2(i-k)-1)!!, \quad \text{(A1)}$$

where $C_i^k$ is the combinatorial number of (16), the operator $(\cdot)!!$ represents double factorial and is given by

$$(2n-1)!! = 1 \times 3 \times 5 \times \ldots \times (2n-1), \quad \text{(A2)}$$

*here we assume* $(-1)!! = 1$.

**Proof.** Mathematical Induction is applied to prove the proposition. Since

$$\sum_{k=1}^{2} C_2^k (2k-3)!!(2(2-k)-1)!! = 3!!, \quad \text{(A3)}$$

$$\sum_{k=1}^{3} C_3^k (2k-3)!!(2(3-k)-1)!! = 5!!, \quad \text{(A4)}$$

(A1) holds when $i=2$ and $i=3$.

Assume that (A1) holds when $i = j-1$, that is,

$$\sum_{k=1}^{j-1} C_{j-1}^k (2k-3)!!(2(j-k)-3)!! = (2j-3)!!. \quad \text{(A5)}$$

Therefore, when $i = j$, we have

$$a = \sum_{k=1}^{j} C_{j-1}^{k-1} (2k-3)!!(2(j-k)-1)!!, \quad \text{(A6)}$$

$$b = \sum_{k=1}^{j} C_{j-1}^{k}(2k-3)!!(2(j-k)-1)!!. \tag{A7}$$

Since

$$C_{j}^{k} = C_{j-1}^{k-1} + C_{j-1}^{k}, (C_{j-1}^{0} = 1), \tag{A8}$$

the following equation is thus given by

$$\sum_{k=1}^{j} C_{j}^{k}(2k-3)!!(2(j-k)-1)!! = a+b. \tag{A9}$$

If $l = k-1$, then we have

$$a = (2j-3)!! + \sum_{l=1}^{j-1} C_{j-1}^{l}(2l-1)!!(2(j-l)-3)!!, \tag{A10}$$

according to (A7), we get the follows

$$b = \sum_{k=1}^{j-1} C_{j-1}^{k}(2k-3)!!(2(j-k)-3)!!(2(j-k)-1), \tag{A11}$$

hence

$$a+b = (2j-3)!! + (2j-2)\sum_{k=1}^{j-1} C_{j-1}^{k}(2k-3)!!(2(j-k)-3)!!. \tag{A12}$$

By considering (A5), we have

$$a+b = (2i-3)!!(2j-1) = (2j-1)!!. \tag{A13}$$

According to (A9), the final equation is given by

$$(2j-1)!! = \sum_{k=1}^{j} C_{j}^{k}(2k-3)!!(2(j-k)-1)!!. \tag{A14}$$

Therefore, (A1) holds. This ends the proof.